\documentclass[12pt,final]{amsart}
\usepackage{amsmath,amssymb}

\usepackage{cite}

\usepackage{color}
\usepackage{soul,graphicx}

\usepackage[color]{showkeys}
\definecolor{refkey}{rgb}{0,0,1}
\definecolor{labelkey}{rgb}{1,0,0}

\begin{document}
\begin{center}
{\bf On the algorithmization of Janashia-Lagvilava matrix spectral factorization method}\\[5mm]
L.Ephremidze$^{1,2}$, F. Saied,$^{1}$ and I. Spitkovsky$^{1}$
\end{center}
\vskip+0.5cm

{\small$^1$ Division of Science and Mathematics, New York University Abu Dhabi (NYUAD),
Saadiyat Island, P.O. Box 129188, Abu Dhabi, United Arab Emirates.\\
$^2$ A. Razmadze Mathematical Institute, I. Javakhishvili Tbilisi State University,
6, Tamarashvili st., Tbilisi 0177, Georgia. E-mail: le23@nyu.edu
           \par  }

\vskip+0.5cm

\small{{\bf Abstract.} We consider three different ways of algorithmization of the Janashia-Lagvilava spectral factorization method. The first algorithm is faster than the second one, however, it is only suitable for matrices of low dimension. The second algorithm, on the other hand, can be applied to matrices of substantially larger  dimension. The third algorithm is a superfast implementation of the method, but only works in the polynomial case under the additional restriction that the zeros of the determinant are not too close to the boundary. All three algorithms fully utilize the advantage of the method which carries out spectral factorization of leading principal submatrices step-by-step. The corresponding results of numerical simulations are reported in order to describe the characteristic features of each algorithm and compare them to other existing algorithms.
\vskip+0.5cm
{\bf Keywords:} Matrix spectral factorization, numerical algorithms.
\vskip+0.5cm

}
\vskip+1.5cm

\section{Introduction}

The Matrix Spectral Factorization (MSF) theorem \cite{Wie57},\cite{HelLow58} asserts
that  if
\begin{equation}\label{1}
S=\begin{pmatrix} s_{11}(t)& s_{12}(t)& \cdots&s_{1r}(t)\\
s_{21}(t)& s_{22}(t)& \cdots&s_{2r}(t)\\
\vdots&\vdots&\vdots&\vdots\\s_{r1}(t)& s_{r2}(t)&
\cdots&s_{rr}(t)\end{pmatrix},
\end{equation}
$|t|=1$, is a  positive definite $($a.e.$)$  matrix function with
integrable entries  defined on the unit circle ${\mathbb T}$ in
the complex plane, $s_{ij}(t)\in L_1({\mathbb T})$, and if the
Paley-Wiener condition
\begin{equation}\label{PW}
\log \det S(t)\in L_1({\mathbb T})
\end{equation}
is satisfied, then $(1)$ admits a spectral factorization
\begin{equation}\label{2}
    S(t)=S^+(t)\big(S^+(t)\big)^*.
\end{equation}
Here the entries of $S^+$ are square integrable functions, $s^+_{ij}\in L_2(\mathbb{T})$, which can be extended analytically inside $\mathbb{T}$, i.e. $s^+_{ij}$ belongs to the Hardy space $H_2$. Furthermore a spectral factor $S^+$ can be selected such that $\det S^+$ is an outer analytic function (see, e.g. \cite{EL10}) and factorization \eqref{2} is unique (up to a constant right unitary multiplier) under these conditions. $S^+$ is unique if we require $S^+(0)$ to be positive definite, and we always assume that it satisfies this condition as well.

In the scalar case, $r=1$, the spectral factor $S^+\in H_2$ can be
explicitly written by the formula
\begin{equation}\label{sc}
 S^+(z)=\exp\left(\frac 1{4\pi}
\int\nolimits_0^{2\pi}\frac{e^{i\theta}+z}{e^{i\theta}-z}\log
S(e^{i\theta})\,d\theta\right).
\end{equation}

If \eqref{1} is a Laurent polynomial matrix
\begin{equation}\label{polm}
S(t)=\sum_{k=-n}^n C_kt^k,\;\; C_k\in{\mathbb C}^{r\times r},
\end{equation}
then the spectral factor
\begin{equation}\label{polm1}
S^+(t)=\sum_{k=0}^n A_kt^k,\;\; A_k\in{\mathbb C}^{r\times r},
\end{equation}
is a polynomial matrix of the same degree $n$ (see e.g. \cite{E} for an elementary proof).

Factorization \eqref{2} was first used in linear prediction theory of multidimensional stationary processes. Nowadays, it is widely known that MSF plays a crucial role in the solution of various applied problems for multiple-input and multiple-output systems in Communications and Control Engineering \cite{Kai99}. Recently MSF became an important step in non-parametric estimations of Granger causality used in Neuroscience \cite{Dhamala},\cite{RoyalA}.
 These applications require the  matrix coefficients of analytic $S^+$ to be determined, at least approximately, for a given matrix function $S$. Therefore, starting with Wiener's original efforts \cite{Wie58} to create a sound computational method of MSF, dozens of different algorithms have appeared in the literature (see the survey papers  \cite{Kuc}, \cite{SayKai} and the references therein, and also \cite{Bott13}, \cite{Jaf} for more recent results).

A novel approach to the solution of the MSF problem, without imposing any additional restriction on $S$ besides the necessary and sufficient condition \eqref{PW} for the existence of spectral factorization, was originally developed by Janashia and Lagvilava in \cite{JL99} for $2\times 2$ matrices. This approach was subsequently  extended to matrices of arbitrary dimension in \cite{IEEE}\footnote{This method obtained USPTO patent recently:\; No. 9,318,232; issued April 19, 2016. }. Results of preliminary numerical simulations based on the proposed method were presented in the same paper \cite{IEEE}. However, a closer look at possible algorithmization ways of this method revealed further advantages. In fact, numerical simulations carried out by the improved algorithms produced much better results than it was reported in \cite{IEEE}. That this development required additional investigations is not surprising, as all methods of MSF are quite demanding and, as it is mentioned in \cite{Kuc}: ``the numerical properties of each method strongly depend on the way it is algorithmized".

In the present paper, after a general description of the Janashia-Lagvilava method (Sections III and IV), we describe three different algorithms of MSF based on this method: JLE-1 (Section VI), JLE-2 (Section VII),  and  JLE-3 (Section VIII). As it was mentioned above, the method is general and also suitable for non-rational matrices. However, since in practical applications the data is finite, we concentrate our attention on the polynomial case. Furthermore, JLE-algorithm 3 is designed only for polynomial matrices \eqref{polm} with the additional restriction that $\det S(t)\not= 0$ for $t\in \mathbb{T}$  (the so-called non-singular case). Its theoretical justification is not yet completed. Nevertheless, due to its superfast speed, we present JLE-3 in the current form. The JLE-algorithm 1 is faster than JLE-2 and it can deal with singular case as well, but it is only suitable for  low dimensional matrices. JLE-algorithm 2 can be applied for much larger matrices, depending on available time and accuracy.  In Section IX, we demonstrate the ability of the method to factorize singular matrices. In Section X, we compare with Wilson's MSF method. The results of provided numerical simulations  are presented in Section XI and concluding remarks are given in Section XII. We emphasize that the proposed MSF method  uses the existing scalar spectral factorization algorithms, whenever they are called for, and does not attempt to improve upon these.

\section{Notation}

Let $\mathbb{D}=\{z\in\mathbb{C}:|z|<1\}$ be the open unit disk, and
$\mathbb{T}=\partial\mathbb{D}$ be the unit circle. As usual,
$L_p=L_p(\mathbb{T})$, $0<p<\infty$, denotes the Lebesgue space
of $p$-integrable complex functions defined on $\mathbb{T}$ ($L_\infty$ is the space of essentially bounded functions). For $p\geq 1$, $\|f\|_p$ is the usual norm. $H_p=H_p(\mathbb{D})$, $0<p\leq\infty$, is the Hardy space of
analytic functions in $\mathbb{D}$ ,
$$
H_p:=\left\{f\in\mathcal{A}(\mathbb{D}):\sup\limits_{r<1}
\int\nolimits_0^{2\pi}|f(re^{i\theta})|^p\,d\theta<\infty\right\}
$$
($H_\infty$ is the space of bounded analytic functions), and
$L_p^+=L_p^+(\mathbb{T})$ denotes the class of their boundary
functions.  A function $f\in H_p$ is called outer, denoted $f\in H_p^O$, if
$$
f(z)=c\cdot  \exp\left(\frac 1{2\pi}
\int\nolimits_0^{2\pi}\frac{e^{i\theta}+z}{e^{i\theta}-z}\log
\big|f(e^{i\theta})\big|\,d\theta\right),\;\;\;\;\;|c|=1.
$$

The $n$th Fourier coefficient of an integrable function $f\in
L_1(\mathbb{T})$ is denoted by $c_k\{f\}$. For $p\geq 1$,
$L_p^+(\mathbb{T})$ coincides with the class of functions from
$L_p(\mathbb{T})$ whose  Fourier coefficients with negative
indices are equal to zero.

The set of trigonometric polynomials is denoted by $\mathcal{P}$,
i.e. $f\in \mathcal{P}$ if $f$ has only a finite number of nonzero
Fourier coefficients. In particular, for integers $m\leq n$,  let $\mathcal{P}_{\{m,n\}}:=\{f\in \mathcal{P}:c_k\{f\}=0
\text{ whenever } k<m \text{ or } k>n\}$ and, for a non-negative integer $N$, let $\mathcal{P}_N^+:=\mathcal{P}_{\{0,N\}}$, $\mathcal{P}_N^-:=\mathcal{P}_{\{-N,0\}}$.
 Obviously, $f\in \mathcal{P}_N^+ \Leftrightarrow \overline{f}\in
\mathcal{P}_N^-$. For a function $f\in L_1$ with Fourier expansion $f\sim\sum_{n\in\mathbb{Z}}c_kt^k$ (or for a formal Fourier series) and positive integer $N$, let $\mathbb{P}_N^+$, $\mathbb{P}_N^-$, and $\mathbb{Q}_N^+$ be the following projection operators:
$$\mathbb{P}_N^+[f]=\sum_{k=0}^N c_kt^k, \mathbb{P}_N^-[f]=\sum_{k=0}^N c_{-k}t^{-k},\text{ and }\mathbb{Q}_N^-[f]=\sum_{k=1}^N c_{-k}t^{-k}$$

If $M $ is a matrix, then $\overline{M}$ denotes the matrix with complex conjugate entries and $M^*:=\overline{M}^T$. Furthermore, $\mathbb{C}^{m\times m}$, $L_p(\mathbb{T})^{m\times m}$, etc., denote the set of $m\times m$ matrices  with the entries from $\mathbb{C}$, $L_p(\mathbb{T})$, etc.
 If $S\in\mathbb{C}^{r\times r}$ is a matrix (function) and $m\leq r$, then  $S_{[m]}$ stands for the upper-left $m\times m$ submatrix of $S$ ($S_{[0]}$ is assumed to be 1) and $S_{[1:\,r,m]}$ stands for $m$th column of $S$. Matrices like $S_{[1:\,r-1,m]}$ or $S_{[1:\,r-1,1:\,m]}$ are defined accordingly. The matrix $S_{]i,j[}$ is obtained from $S$ by deleting the $i$th row and $j$th column.

 For a polynomial $p(t)=\sum_{k=0}^k p_kt^k$, let $\|p\|=\sup_{0\leq k\leq n}|p_k|$, and for a polynomial matrix $P=\big(P_{ij}\big)_{i,j=1}^r$, let $\|P\|=\sup_{1\leq i,j\leq r}\|P_{ij}\|$.

 A matrix $M\in\mathbb{C}^{r\times r}$ is called positive definite if $X^*MX>0$ for all $0\not=X\in \mathbb{C}^{r\times1}$, and $S\in L_1(\mathbb{T})^{r\times r}$ is called positive definite if it is positive definite for a.a. $t\in \mathbb{T}$.

 A matrix function $U\in L^\infty(\mathbb{T})^{r\times r}$ is called unitary if
 \begin{equation}
\label{Ud}
 U(t)U^*(t)=I_r \;\;\text{a.e.},
\end{equation}
where $I_r$ stands for $r\times r$ identity matrix.

$\mathbf{0}_{r\times m}$ and $\mathbf{1}_{r\times m}$ stand for $r\times m$ matrices with all entries equal to $0$ and $1$, respectively. Using Matlab's notation, if $A\in \mathbb{C}^{r\times m_1}$ and $B\in \mathbb{C}^{r\times m_2}$, then $[A\;B]$ is $r\times(m_1+m_2)$ matrix, while if $A\in \mathbb{C}^{r_1\times m}$ and $B\in \mathbb{C}^{r_2\times m}$, then $[A\,;\;B]=[A^T\;B^T]^T$ is $(r_1+r_2)\times m$ matrix.

For a column vector $\mathbf{a}=[a_0\,a_1\,\cdots\,a_{l}]^T\in \mathbb{C}^{(l+1)\times 1}$ and a positive integer $m\in\mathbb{N}$, let $T(\mathbf{a};m)$ be the $(l+m+1)\times (m+1)$ Toeplitz  matrix with the first column $[\mathbf{a}\,; \mathbf{0}_{m\times 1}]\in\mathbb{C}^{l+m+1}$ and the first row $[a_0\,\mathbf{0}_{1\times m}]\in\mathbb{C}^{1\times (m+1)}$.

We say that a sequence of matrix functions $S_n$, $n=1,2,\ldots$ is convergent to a matrix function $S$ (in some sense) if the entries of $S_n$ are convergent to the corresponding entries of $S$ (in this sense).

Finally, $\delta_{ij}$ stands for the  Kronecker delta, i.e. $\delta_{ij}=1$ if $i=j$ and $\delta_{ij}=0$ otherwise.

\section{General description of the method }

The first step of the MSF method proposed in \cite{IEEE}  is the triangular factorization of \eqref{1}
\begin{equation}
\label{S1}
S(t)=M(t)M^*(t),
\end{equation}
where $M(t)$ is the lower triangular matrix
\begin{equation}
\label{S2}
M(t)=\begin{pmatrix}f^+_1(t)&0&\cdots&0&0\\
        \xi_{21}(t)&f^+_2(t)&\cdots&0&0\\
        \vdots&\vdots&\vdots&\vdots&\vdots\\
        \xi_{r-1,1}(t)&\xi_{r-1,2}(t)&\cdots&f^+_{r-1}(t)&0\\
        \xi_{r1}(t)&\xi_{r2}(t)&\cdots&\xi_{r,r-1}(t)&f^+_r(t)
        \end{pmatrix},
\end{equation}
$\xi_{ij}\in L_2(\mathbb{T})$, $f_i^+\in H_2^O$. The spectral factor $S^+$ is represented in the form
\begin{equation}\label{S3}
S^+(t)=M(t)\mathbf{U}_2(t)\mathbf{U}_3(t)\ldots\mathbf{U}_r(t)\cdot U.
\end{equation}
Here each $\mathbf{U}_m$ is a block matrix function
\begin{equation}\label{UB}
\mathbf{U}_m(t)=\begin{pmatrix}U_{m}(t)&\mathbf{0}_{m\times(r-m)}\\\mathbf{0}_{(r-m)\times m}&I_{r-m}\end{pmatrix},
\end{equation}
where $U_m(t)$ is a special unitary matrix function of the form
\begin{equation}\label{S4}
U_m(t)=\begin{pmatrix}u_{11}(t)&u_{12}(t)&\cdots&u_{1,m-1}(t)&u_{1m}(t)\\
                 u_{21}(t)&u_{22}(t)&\cdots&u_{2,m-1}(t)&u_{2m}(t)\\
           \vdots&\vdots&\vdots&\vdots&\vdots\\
           u_{m-1,1}(t)&u_{m-1,2}(t)&\cdots&u_{m-1,m-1}(t)&u_{m-1,m}(t)\\[3mm]
           \overline{u_{m1}(t)}&\overline{u_{m2}(t)}&\cdots&\overline{u_{m,m-1}(t)}&\overline{u_{mm}(t)}\\
           \end{pmatrix},
\end{equation}
 with
\begin{equation} \label{S45}
u_{ij}\in L^\infty_+,\;\; \text{ and }\;\; \det U(t)=1 \text{ a.e.}
\end{equation}
(for reasons explained in \cite{EL2014} such matrices can as well be called ``wavelet matrices").  Furthermore, for each $m=2,3,\ldots,r$,
\begin{equation}\label{75}
S_{[m]}^+=\big(M\mathbf{U}_2\mathbf{U}_3\ldots\mathbf{U}_m\big)_{[m]}
\end{equation}
is a spectral factor of $S_{[m]}$. In particular, $S^+_0:=M\mathbf{U}_2\mathbf{U}_3\ldots\mathbf{U}_r$ is a spectral factor of \eqref{1}, and the constant unitary matrix $U$ in \eqref{S3} makes $S^+$ positive definite in the origin, namely (see \cite[formula (54)]{EJL11})
\begin{equation}\label{Uc}
U=\big(S^+_0(0)\big)^{-1}\sqrt{S^+_0(0)(S^+_0(0))^*}.
\end{equation}

To obtain unitary matrix function \eqref{S4} for each $m=2,3,\ldots,r$ recurrently, we  consider a matrix function
\begin{equation}
\label{S5}
F_m(t)=\begin{pmatrix}1&0&0&\cdots&0&0\\
          0&1&0&\cdots&0&0\\
           0&0&1&\cdots&0&0\\
           \vdots&\vdots&\vdots&\vdots&\vdots&\vdots\\
           0&0&0&\cdots&1&0\\
           \zeta_{1}(t)&\zeta_{2}(t)&\zeta_{3}(t)&\cdots&\zeta_{m-1}(t)&f^+_m(t)
           \end{pmatrix},
\end{equation}
where the last row of \eqref{S5} consists of the first $m$ entries of the $m$th row of the product
\begin{equation}\label{121}
M_{m-1}:=M\mathbf{U}_2\mathbf{U}_3\ldots\mathbf{U}_{m-1},
\end{equation}
and then obtain a matrix function \eqref{S4}, \eqref{S45} such that (see \cite[Lemma 4]{EJL11})
\begin{equation}\label{16.5}
F_mU_m\in L_2^+(\mathbb{T})^{m\times m}.
\end{equation}
Particularly, we have
\begin{equation}\label{17.1}
\big(M_{m-1}\big)_{[m]}=\left[\begin{matrix}&   S_{[m-1]}^+(t)& & \begin{matrix}0\\0\\ \vdots\\0\end{matrix}\\
\zeta_1(t) &  \ldots& \zeta_{m-1}(t)&  f^+_m(t)
\end{matrix}\right]=\left[\begin{matrix}&   S_{[m-1]}^+(t)& & \begin{matrix}0\\0\\ \vdots\\0\end{matrix}\\
0 &  \ldots& 0&  1
\end{matrix}\right] F_m(t)
\end{equation}
and
\begin{equation}\label{17.2}
S_{[m]}^+(t)=\left[\begin{matrix}& &  S_{[m-1]}^+(t)& & \begin{matrix}0\\0\\ \vdots\\0\end{matrix}\\
0 & 0 & \ldots& 0&  1
\end{matrix}\right] F_m(t) U_m(t).
\end{equation}

In order to achieve \eqref{16.5}, one needs to consider the following system of conditions (see \cite[formula (15)]{IEEE})
\begin{equation}\label{I15}
\begin{cases} \zeta_1(t)x^+_m(t)-f_m^+(t)\overline{x^+_1(t)}\in L_2^+,\\
              \zeta_2(t)x^+_m(t)-f_m^+(t)\overline{x^+_2(t)}\in L_2^+,\\
              \vdots\\
              \zeta_{m-1}(t)x^+_m(t)-f_m^+(t)\overline{x^+_{m-1}(t)}\in L_2^+,\\
              \zeta_1(t)x^+_1(t)+\zeta_2(t)x^+_2(t)+\ldots+\zeta_{m-1}(t)x^+_{m-1}(t)
              +f_m^+(t)\overline{x^+_m(t)}\in L_2^+,
         \end{cases}
\end{equation}
and columns of \eqref{S4} are $m$ independent solutions of \eqref{I15}.

To construct \eqref{S4} approximately the following procedures should be performed:

For a large positive $N$, let $F_m^{\{N\}}$ be the matrix function \eqref{S5} with the last row replaced by
$$
(\zeta_1^{\{N\}},\zeta_2^{\{N\}},\ldots,\zeta_{m-1}^{\{N\}},f_m^+),
$$
where
$$
\zeta_j^{\{N\}}(t):=\sum_{k=-N}^\infty c_k\{\zeta_j\}t^{k},\;\;\;j=1,2,\ldots,m-1.
$$
Then one can find the unitary matrix function $U_m^{\{N\}}$ of the form \eqref{S4} such that $\det U_m^{\{N\}}(t)=1$, $U_m^{\{N\}}(1)=I_m$, $u_{ij}\in\mathcal{P}_N^+$ and $F_m^{\{N\}}U_m^{\{N\}}\in \mathcal{P}_N^+$ (see \cite[Theorem 1]{IEEE}). In particular, the columns of $U_m^{\{N\}}$ are $m$ independent solutions of the  system \eqref{I15} where $\zeta_1,\zeta_2,\ldots,\zeta_{m-1}$    are replaced by $\zeta_1^{\{N\}},\zeta_2^{\{N\}},\ldots,\zeta_{m-1}^{\{N\}}$, and they can be actually found by solving a single system of $(N+1)\times (N+1)$ linear algebraic equations with $m$ different right-hand sides (see the proof of Theorem 1 in \cite{IEEE}). Details of the computation are given in Section IV.

One can prove that $U_m^{\{N\}}\to U_m$ at least in measure as $N\to\infty$, which guarantees that (see \cite[Theorem 2]{EJL11})
$$
M\mathbf{U}_2\mathbf{U}_3\ldots\mathbf{U}_{m}^{\{N\}}\to M\mathbf{U}_2\mathbf{U}_3\ldots\mathbf{U}_{m}\;\;\text{ in } L_2.
$$

\section{Construction of wavelet matrices}

In this section we provide the details of computation of the unitary matrix function $U_N:=U_m^{\{N\}}$ for a given matrix function \eqref{S5}. $N$ and $m$ are assumed fixed throughout this section.

Let
$$\mathbb{P}_N^+[f_m^+](t)=\sum_{k=0}^N d_kt^k,\;\;\;\mathbb{Q}_N^-[\zeta_i](t)=\sum_{k=1}^N\gamma_{in}t^{-k},\;\;\text{ and } \;\;\mathbb{P}_N^+[1/f_m^+](t)=\sum_{k=0}^Nb_kt^k.$$
(Note that the knowledge of $\mathbb{P}_N^+[f_m^+]$ is sufficient to determine $\mathbb{P}_N^+[1/f_m^+]$.) Suppose $D^{-1}$ is the upper triangular Toeplitz matrix with the first row
\begin{equation}\label{data1}
(b_0,b_1,\ldots,b_N),
\end{equation}
and $\Gamma_i$, $i=1,2,\ldots,m-1$ is the upper triangular Hankel matrix withe the first row
\begin{equation}\label{data2}
(0,\gamma_{i,1},\gamma_{i,2},\ldots,\gamma_{iN})
\end{equation}
(see \cite[(26)]{IEEE}) and let
\begin{equation}\label{tht}
\Theta_i=D^{-1}\,\Gamma_i\,,\;\;i=1,2,\ldots,m-1.
\end{equation}
Note that $\Theta_i$ is the upper triangular Hankel matrix (see \cite[(33)]{IEEE}) with the first row
\begin{equation}\label{Lb}
\Lambda_i:=(\eta_{i0}, \eta_{i1},\ldots,\eta_{iN}),
\end{equation}
where $\sum_{k=0}^N \eta_{in}t^{-k}=\mathbb{P}_N^-\left[\sum_{k=0}^Nb_kt^k\cdot\sum_{k=1}^N\gamma_{in}t^{-k}\right]$.

Take
\begin{equation}\label{dt}
\Delta=\sum_{i=1}^{m-1}\Theta_i\Theta_i^*+I_{N+1},
\end{equation}
which is a positive definite matrix (with all eigenvalues $\geq 1$), and solve the same system of equations (see \eqref{Lb})
\begin{equation}\label{eqdet}
\Delta X=\Lambda_i^T
\end{equation}
with $m$ different right hand sides corresponding to $i=1,2,\ldots,m$. Here it is assumed that $\Lambda_m=(1,0,0,\ldots,0)$. The matrix \eqref{dt} has a displacement structure of rank $m$, namely
$$\Delta-Z\Delta Z^*=\sum_{i=1}^{m-1}\Lambda_i\Lambda_i^*+\mathcal{E}\mathcal{E}^*$$
 has rank $m$, where $Z$ is the upper triangular $(N+1)\times(N+1)$
matrix with 1's on the first up-diagonal and 0's elsewhere (i.e. a
Jordan block with eigenvalue $0$) and ${\mathcal E}=(0,0,\ldots,0,1)^T \in \mathbb{C}^{N+1,1}$ (see \cite[Appendix]{IEEE}). Therefore its triangular factorization $\Delta=LDL^*$ can be achieved in $O(mN^2)$ operations instead of $O(N^3)$ as explained e.g. in \cite[Appendix F]{Kai99} without even constructing the matrix $\Delta$ (just using the $(N+1)\times m$ matrix $[\Lambda_1,\Lambda_2,\ldots,\Lambda_{m-1},\mathcal{E}]$).

Let the solution of \eqref{eqdet} be $X_i=(a_{i0},a_{i1},\ldots, a_{iN})^T$, and denote
$$
v_{mi}(t):=\sum_{k=0}^N a_{in}t^k,\;\;\;\;i=1,2,\ldots,m,
$$
Suppose also
$$
v_{ij}(t)=\mathbb{P}_N^+\left[\sum_{k=0}^N \overline{\eta}_{n}t^k\cdot\sum_{k=0}^N \overline{a_{in}}t^{-k}\right]-\delta_{ij},
$$
$1\leq i\leq m-1$, $1\leq j\leq m$, and let
$$
V(t)=\begin{pmatrix}v_{11}(t)&v_{12}(t)&\cdots&v_{1,m-1}(t)&v_{1m}(t)\\
                 v_{21}(t)&v_{22}(t)&\cdots&v_{2,m-1}(t)&v_{2m}(t)\\
           \vdots&\vdots&\vdots&\vdots&\vdots\\
           v_{m-1,1}(t)&v_{m-1,2}(t)&\cdots&v_{m-1,m-1}(t)&v_{m-1,m}(t)\\[3mm]
           \overline{v_{m1}(t)}&\overline{v_{m2}(t)}&\cdots&\overline{v_{m,m-1}(t)}&\overline{v_{mm}(t)}\\
           \end{pmatrix}.
$$
Then  (see \cite[(51)]{IEEE})
$$
U_N(t)=V(t)\cdot V^{-1}(0).
$$
It is proved in \cite{IEEE} that $V(0)$ is nonsingular and the condition number of this matrix is estimated in \cite{ESS}.

\section{A shortcut in the recursive  step}
As it was mentioned in Section III, in order to perform $m$th recursive  step in the proposed MSF method, we need only to consider
$$
S_{[m-1]}^+=\big(M_{m-1}\big)_{[m-1]}
$$
(see \eqref{121} and \eqref{75}), which has already been constructed (at least approximately) and the first $m$ entries in the $m$th row of $M_{m-1}$
\begin{equation}\label{tr1}
\big(M_{m-1}\big)_{[m,1:m]}=(\zeta_1,\zeta_2,\ldots,\zeta_{m-1},f^+_m)
\end{equation}
(see \eqref{S5}). Because of the block structure of matrices in \eqref{UB}, the entry $f^+_m$ is the same as in \eqref{S2}. Thus it can be computed by the formula (see \cite[formula (56)]{IEEE})
\begin{equation}\label{fm}
f^+_m=\frac{\big(\det S_{[m]}\big)^+}{\big(\det S_{[m-1]}\big)^+}
\end{equation}
(\,$(\cdot)^+$ stands for the scalar spectral factorization \eqref{sc}\,).

Since $S=M_{m-1}M_{m-1}^*$ (see \eqref{S1}, \eqref{121}, and \eqref{Ud}) and particularly
\begin{equation}\label{tr2}
(S)_{[m]}=\big(M_{m-1}\big)_{[m]}\big(M_{m-1}\big)_{[m]}^*
\end{equation}
(see \eqref{17.1}), we have
\begin{equation}\label{Cr1}
S_{[m-1]}^+\cdot\big(\zeta_1,\zeta_2,\ldots,\zeta_{m-1}\big)^*=S_{[1:\,m-1,\,m]}.
\end{equation}
Therefore, instead of computing matrices $M_m$ for each $m=2,3,\ldots,r-1$ by \eqref{121}, we can directly compute the  entries $\zeta_1,\zeta_2,\ldots,\zeta_{m-1}$ from \eqref{Cr1}.

Having computed the functions $\zeta_1,\zeta_2,\ldots,\zeta_{m-1}$, one can find $|f^+_m|^2$ from the formula (see \eqref{tr2})
\begin{equation}\label{275}
\sum_{j=1}^{m-1}|\zeta_j|^2+|f_m^+|^2=s_{mm}.
\end{equation}
Therefore, an alternative way of computing \eqref{fm} is the scalar spectral factorization of $s_{mm}-\sum_{j=1}^{m-1}|\zeta_j|^2$.

In the next three sections we present three different implementations of the described algorithm for polynomial data \eqref{polm}, followed by the results of corresponding numerical simulations.

\section{JLE-algorithm 1}

This algorithm relies on computation of polynomial matrix determinant. Namely, for a polynomial matrix of order $n$
\begin{equation}\label{pol2}
P(t)=\sum_{k=0}^n B_kt^k,\;\; B_k\in{\mathbb C}^{m\times m},
\end{equation}
$\det P$ is a polynomial of order $mn$. Therefore,  having evaluated $\det P(t)$  at $mn+1$ DFT nodes $t_l=\exp\left(\frac{2\pi il}{mn+1}\right)$, $l=0,1,\ldots,mn$,  the coefficients of $\det P$ can be computed by interpolation, namely computing the inverse DFT of $[\det P(t_0), \ldots,\det P(t_{mn})]$.

This algorithm of polynomial matrix determinant computation is fast and accurate for matrices of small dimension. However,  the algorithm suffers from severe round-off errors and the accuracy is destroyed for large dimensional matrices.  For example, with a standard double precision in Matlab, we have found a computation error in the formula
$$\| \det(P_1P_2)-\det P_1\,\det P_2\|$$
as small as $10^{-8}$ for randomly selected polynomial matrices $P_1$ and $P_2$ of degree $n=10$ and dimension $m=10$, and as large as $10^9$ for ones with $n=20$ and $m=15$. The reason of such increase is that the coefficients of $\det P$ become very large (at least for randomly selected coefficients $B_k$ in \eqref{pol2}) and floating point machine arithmetic  loses significant digits. Therefore JLE-algorithm 1 (with input \eqref{polm} and output \eqref{polm1}) is suitable for small dimensional matrices  ($r<20$ and $n<25$).  Its basic computational procedures are described below.

\smallskip
{\bf Procedure 1.} Compute the diagonal entries of the triangular factor \eqref{S2} by the formula \eqref{fm}, where $m=1,2,\ldots,r$. Each $f_m^+$ can be represented as a rational function $p_m/q_m$, where $p_m\in\mathcal{P}_{mn}^+$ and $q_m\in\mathcal{P}_{(m-1)n}^+$. In addition, the denominator is free of zeros inside $\mathbb{T}$, and $f_m^+$ is free of poles on $\mathbb{T}$ (since $f_m^+\in L_2^+(\mathbb{T})$).

For the scalar spectral factorization of $\det S_{[m]}$, we first apply exp-log implementation by using FFT \cite{SFLP} and then we improve the accuracy by using 4-5 iterations of Wilson's scalar factorization algorithm \cite{Wil1}.

\smallskip
{\bf Procedure 2.} For $m=2,3,\ldots,r$, assume that $S_{[m-1]}^+$ has already been (approximately) constructed as an $(m-1)\times(m-1)$ polynomial matrix of degree $n$ and perform the following steps.

{\bf Step 1.} Compute $\zeta_j$, $j=1,2,\ldots,m-1$, by the Cramer's rule from equation \eqref{Cr1}.  In particular, each $\overline{\zeta_j}$ will be of the form $p/q$, where $p\in\mathcal{P}_{\{-n,\,(m-1)n\}}$ and $q\in\mathcal{P}_{(m-1)n}^+$, again with $q$ free of zeros inside $\mathbb{T}$ and $\zeta_j$ free of poles on $\mathbb{T}$. Note that $\zeta_j$-s can be computed in parallel.

{\bf Step 2.} Select a large positive integer $N$. Theoretically, as $N\to\infty$, the computed spectral factor $\hat{S}_{[m]}^+$ converges to exact ${S}_{[m]}^+$ (assuming that all previous factors including $S_{[m-1]}^+$ are computed exactly). However, in practise we never achieve an exact result. Nevertheless, the accuracy
\begin{equation}\label{iac}
\|{S}_{[m]}-\hat{S}_{[m]}^+\big(\hat{S}_{[m]}^+\big)^*\|
\end{equation}
can be controlled and the value of $N$ can be increased, if necessary, at each intermediate stage, in order to achieve a satisfactory approximation in the  final result.

{\bf Step 3.} From obtained representations of $\zeta_j$, $j=1,2,\ldots,m-1$, and $f_m^+$ as rational functions, find
$$
\zeta_j^{\{N\}}:=\mathbb{Q}_N^-[\zeta_j]+\mathbb{P}_n^+[\zeta_j]=\sum_{k=-N}^n c_k\{\zeta_j\}t^{k}
$$
and
$$
f_m^{\{N\}}:=\mathbb{P}_{N+n}^+[f_m^+]=\sum_{k=0}^{N+n}c_k\{f_m^+\}t^k.
$$
We do this by the standard division algorithm  of two polynomials, utilizing the advantages of denominator being  free from zeros inside $\mathbb{T}$ and function having no poles on $\mathbb{T}$.

{\bf Step 4.} Using $(\zeta_1^{\{N\}},\zeta_2^{\{N\}},\ldots,\zeta_{m-1}^{\{N\}}, f_m^{\{N\}})$ as the last row of \eqref{S5}, construct a unitary matrix function $U_N:=U_m^{\{N\}}$ as it is described in Section IV.

{\bf Step 5.} Consider the product
$$
S_{[m]}^+\approx \begin{pmatrix}&&&&0\\&&S_{[m-1]}^+&&\vdots\\&&&&0\\
                {\zeta}_{1}^{\{N\}}& {\zeta}_{2}^{\{N\}}&\ldots& {\zeta}_{m-1}^{\{N\}}& {f}^{\{N\}}_m
           \end{pmatrix}
           \begin{pmatrix}u_{11}&u_{12}&\cdots&u_{1m}\\
                  \vdots&\vdots&\vdots&\vdots\\
           u_{m-1,1}&u_{m-1,2}&\cdots&u_{m-1,m}\\[3mm]
           \overline{u_{m1}}&\overline{u_{m2}}&\cdots&\overline{u_{mm}}\\
           \end{pmatrix}
$$
(the last matrix is $U_m^{\{N\}}$), where all coefficients of polynomials in the right-hand side product with indices outside the range $[0,n]$ are neglected (since we know that the exact $S_{[m]}^+$ is matrix polynomial of degree $n$). Therefore, $S_{[m-1]}^+$ can be separately multiplied by the first $m-1$ rows of $U_m^{\{N\}}$ and then its last row can be multiplied by $U_m^{\{N\}}$.

\smallskip
{\bf Procedure 3.} For $m=r$, $S_{[r]}^+$ is an approximate spectral factor of $S$. We can multiply $S_{[r]}^+$  by the constant unitary matrix $U$ defined by \eqref{Uc} (taking $S_{[r]}^+$ instead of $S_{0}^+$ ) to obtain $S^+$.

\section{JLE-algorithm 2}

In this implementation, computations of polynomial matrix determinants are avoided. Consequently much higher dimensional matrices can be factorized accurately by this algorithm at the expense of large computer memory usage.

{\bf Procedure 1.} Compute a scalar spectral factor $f_1^+$ of $s_{11}$ by using the same exp-log and Wilson's methods as in Procedure 1 of JLE-algorithm 1.

\smallskip
{\bf Procedure 2.} For $m=2,3,\ldots,r$, assume that $S_{[m-1]}^+$ has already been (approximately) constructed as an $(m-1)\times(m-1)$ polynomial matrix of degree $n$ and perform the following steps.

{\bf Step 1.} Take a large number of DFT nodes, usually $2^\kappa$, where $10\leq\kappa\leq 23$: $t_l=\exp\left(\frac{2\pi il}{2^\kappa}\right)$, $l=0,1,\ldots,2^\kappa-1$. This $\kappa$ becomes another tuning parameter in the algorithm (along with $N$), which can be selected and changed during recursive  steps in order to improve the accuracy \eqref{iac}.

{\bf Step 2.} For each node $t_l$, $l=0,1,\ldots,2^\kappa-1$, evaluate the matrices $S_{[m-1]}^+(t_l)$ and $S_{[1:\,m-1,\,m]}(t_l)$, and solve the following system of linear equations (see \eqref{Cr1}):
\begin{equation}\label{sys2}
S_{[m-1]}^+(t_l)\cdot X=S_{[1:\,m-1,\,m]}(t_l).
\end{equation}
We have $\big(\zeta_1(t_l),\zeta_2(t_l),\ldots,\zeta_{m-1}(t_l),\big)=X_l^*$, where $X_l$ is the solution of \eqref{sys2}.

If it happens that the system \eqref{sys2} is singular or ill conditioned, then we can apply the continuity of functions $\zeta_j$ and assume that $X_l=X_{l-1}$.

When standard routines are well optimized (as it is in Matlab), this step is not as time-consuming as it might appear at the first glance.

{\bf Step 3.} Compute $|f_m^+(t_l)|^2$, $l=0,1,\ldots,2^\kappa-1$, from the formula \eqref{275}

{\bf Step 4.} Select a large positive integer $N$,  and using the values of $|f_m^+|^2$ at DFT nodes, perform an approximate scalar spectral factorization to reconstruct
$$
f_m^{\{N\}}:=\sum_{k=0}^{N+n}c_k\{f_m^+\}t^k.
$$
For this step, one can use the exp-log method of scalar spectral factorization which utilizes the boundary values of a spectral density.

The integer $N$ has a natural bound  $2^\kappa-n$ in this situation, however an optimal ratio (from 1/10 to 1/50) of $N/2^\kappa$ should be selected in order to achieve a good accuracy.

{\bf Step 5.}
From the values of $\zeta_j$ at DFT nodes $t_l$, $l=0,1,\ldots,2^\kappa-1$, reconstruct (approximately)
\begin{equation}\label{zeta1}
\zeta_j^{\{N\}}:=\sum_{k=-N}^n c_k\{\zeta_j\}t^{k}
\end{equation}
by using the inverse FFT and selecting corresponding coefficients.

\smallskip

The remaining steps are the same as Steps 4 and 5 in JLE-algorithm 1, including Procedure 3.

\section{JLE-algorithm 3}
This implementation utilizes formulas \eqref{17.2}, \eqref{fm}, and \eqref{Cr1} for $m=r$:

\begin{equation} \label{dr7}
S^+(t)=\left[\begin{matrix}& &  S_{[r-1]}^+(t)& & \begin{matrix}0\\0\\ \vdots\\0\end{matrix}\\
\zeta_1(t) & \zeta_2(t) & \ldots& \zeta_{r-1}(t)&  f_r^+(t)
\end{matrix}\right]
\left[\begin{matrix}u_{11}(t)&u_{12}(t)&\cdots&u_{1r}(t)\\
                 u_{21}(t)&u_{22}(t)&\cdots&u_{2r}(t)\\
           \vdots&\vdots&\vdots&\vdots\\
           u_{r-1,1}(t)&u_{r-1,2}(t)&\cdots&u_{r-1,r}(t)\\[3mm]
           \overline{u_{r1}(t)}&\overline{u_{r2}(t)}&\cdots&\overline{u_{rr}(t)}\\
           \end{matrix}\right],
\end{equation}
\begin{equation} \label{dr8}
 f_r^+(t)=\det S^+(t)/\det S^+_{[r-1]}(t),
\end{equation}
and
\begin{equation} \label{dr9}
[\zeta_1(t), \zeta_2(t), \ldots, \zeta_{r-1}(t)]\cdot\big(S^+_{[r-1]}(t)\big)^*=S_{[r,1:\,r-1]}(t).
\end{equation}

Let $U(t)=U_r(t)$ be the last matrix in \eqref{dr7}.
Then,  for $j\leq r$, it follows from \eqref{dr7} that
\begin{equation} \label{dr10}
 S^+_{[r-1]}(t)\cdot U_{[1:\,r-1,\,j]}(t)=S^+_{[1:\,r-1,\,j]}(t),
\end{equation}
and furthermore
\begin{equation} \label{dr11}
 S^+_{[r-1]}(t)\cdot U_{]r,j[}(t)=S^+_{]r,j[}(t).
\end{equation}

Since $U(t)$ is a unitary matrix ($U^{-1}(t)=U^*(t)$) and $\det U(t)=1$, it follows that $u_{r,j}(t)=\det U_{]r,j[}(t)$ and, taking into account \eqref{dr11}, we get
\begin{equation} \label{dr12}
 u_{r,j}(t)=\det S^+_{]r,j[}(t)/\det S^+_{[r-1]}(t).
\end{equation}

It also follows from \eqref{dr7} that
\begin{equation} \label{dr13}
[\zeta_1(t), \zeta_2(t), \ldots, \zeta_{r-1}(t)]\cdot U_{[1:\,r-1,\,j]}(t)+f_r^+(t)\overline{u_{r,j}(t)}=S^+_{rj}(t).
\end{equation}

Substituting into \eqref{dr13} $[\zeta_1,  \ldots, \zeta_{r-1}]=S_{[r,1:\,r-1]}\cdot\big(S^+_{[r-1]}\big)^{-*}$ (see \eqref{dr9}\,), $  U_{[1:\,r-1,\,j]}=\big(S^+_{[r-1]}\big)^{-1}\cdot S^+_{[1:\,r-1,\,j]}$ (see \eqref{dr10}\,),  \eqref{dr8}, and \eqref{dr12},  and taking into account that $S_{[r-1]}=S^+_{[r-1]}\big(S^+_{[r-1]}\big)^*$, we get

$$
S_{[r,1:\,r-1]}(t)\cdot\big(S_{[r-1]}(t)\big)^{-1}\cdot S^+_{[1:\,r-1,\,j]}(t)+\frac{\det S^+(t)\overline{\det S^+_{]r,j[}(t)}}{\det S_{[r-1]}(t)}=
S^+_{rj}(t).
$$
Consequently,
\begin{equation} \label{dr16}
S_{[r,1:\,r-1]}(t)\cdot\mathop{\rm Cof}\big\{S_{[r-1]}(t)\big\}^{T}\cdot S^+_{[1:\,r-1,\,j]}(t)+{\big({\det S(t)\big)^+}\cdot\overline{\det S^+_{]r,j[}(t)}}=
S^+_{rj}(t){\det S_{[r-1]}(t)},
\end{equation}
where it is assumed that $\big(\det S(t)\big)^+$ can be found from  ${\det S(t)}$, as the problem is reduced to the scalar spectral factorization.

 In the equation \eqref{dr16},  $S_{[r,1:\,r-1]}$, $\mathop{\rm Cof}\big\{S_{[r-1]}\big\}^{T}$, $\big(\det S(t)\big)^+$ and ${\det S_{[r-1]}}$ are  assumed to be the known (matrix) functions, and $S^+_{[1:\,r-1,\,j]}$, ${\det S^+_{]r,j[}}$, and $S^+_{rj}$ are  unknown (matrix) functions.

 Assume now that $S$ is a matrix polynomial of degree $n$ (see \eqref{polm}), i.e. $S\in(\mathcal{P}_{\{-n,n\}})^{r\times r}$.
Let us observe that for functions in \eqref{dr16} we have:
\begin{gather*}
{S}_{[r,1:\,r-1]}\in(\mathcal{P}_{\{-n,n\}})^{1\times(r-1)}; {\mathop{\rm Cof}\big\{S_{[r-1]}\big\}}^{T}\in (\mathcal{P}_{\{-n(r-2), n(r-2)\}})^{(r-1)\times(r-1)};\\ {S^+_{[1:\,r-1,\,j]}}\in(\mathcal{P}_{\{0,n\}})^{(r-1)\times 1}; {(\det S)^+}\in \mathcal{P}_{\{0,rn\}};  {\overline{\det S^+_{]r,j[}}}\in\mathcal{P}_{\{-(r-1)n,0\}}; {S^+_{rj}}\in\mathcal{P}_{\{0,n\}},
 \end{gather*}
and ${\det S_{[r-1]}}\in\mathcal{P}_{\{-(r-1)n, (r-1)n\}}$.  Thus all products in \eqref{dr16} have the range of indices of (nonzero) Fourier coefficients in $[-(r-1)n,rn]$. If we equate the corresponding coefficients in these products, we get $2rn-n+1$ linear algebraic equations with respect to coefficients of unknown (matrix) polynomials $S^+_{[1:\,r-1,\,j]}$, ${\det S^+_{]r,j[}}$, and $S^+_{rj}$. The total number of these coefficients is $(r-1)(n+1)+\{(r-1)n+1\}+(n+1)=2rn-n+r+1$.

We can factorize $S(t)$ at a single point on the unit circle, say $t=1$, and getting the representation $S(1)=S^+(1)\big(S^+(1)\big)^*$, we can assume that $[S^+_{[1:\,r-1,\,j]}(1)\,S^+_{r,j}(1)]^T$ is the $j$-th column of $S^+(1)$. This gives the additional $r$ conditions on coefficients of (matrix) polynomials $S^+_{[1:\,r-1,\,j]}$  and $S^+_{r,j}$, and thus additional $r$ equations. In the end we get the same number of linear equations and unknowns $2rn-n+r+1$.

The basic computational procedures of the algorithm are described below.

\smallskip

{\bf Step 1.} Compute the polynomial determinants $\det S(t)$ and $\det S_{[r-1]}(t)$ by the method described in JLE-1.

\smallskip

{\bf Step 2.} Compute the scalar spectral factor $\big(\det S(t)\big)^+$ by the method described in Procedure 1 of JLE-1.

\smallskip

{\bf Step 3.} Compute $\mathop{\rm Cof}\big\{S_{[r-1]}(t)\big\}^{T}$ by evaluating it at $N=2n(r-1)+1$ DFT nodes $t_l=\exp\left(\frac{2\pi il}{N}\right)$, $l=0,1,2,\ldots,N-1$, by the formula $\mathop{\rm Cof}\big\{S_{[r-1]}(t_l)\big\}^{T}=\det S_{[r-1]}(t_l)\big(S_{[r-1]}(t_l)\big)^{-1}$ and then use the inverse Fourier transform.

\smallskip

{\bf Step 4.} Multiply matrix polynomials $S_{[r,1:\,r-1]}$ and $\mathop{\rm Cof}\big\{S_{[r-1]}\big\}^T$.

\smallskip

Let $\big(\det S(t)\big)^+=\sum_{k=0}^{rn}a_kt^k$,\;\;\; $t^{(r-1)n}\det S_{[r-1]}(t)=\sum_{k=0}^{2(r-1)n}b_kt^k$, and
$$
t^{(r-1)n}S_{[r,1:\,r-1]}(t)\mathop{\rm Cof}\big\{S_{[r-1]}\big\}^T(t)=\sum_{k=0}^{2(r-1)n}C_kt^k=\big[\sum_{k=0}^{2(r-1)n}c_k^{\{1\}}t^k
\cdots \sum_{k=0}^{2(r-1)n}c_k^{\{r-1\}}t^k\big],
$$
$C_k\in\mathbb{C}^{1\times(r-1)}$, $c_k^{\{j\}}\in\mathbb{C}$. Introduce also the notation: $\mathbf{a}=[a_0\,a_1\,\cdots\,a_{rn}]^T\in \mathbb{C}^{(2rn+1)\times 1}$; $\mathbf{b}=[b_0\,b_1\,\cdots
b_{2(r-1)n}]^T\in  \mathbb{C}^{(2(r-1)n+1)\times 1}$; $\mathbf{c}^{\{j\}}=[c_0^{\{j\}}\,c_1^{\{j\}}\,\cdots c_{2(r-1)n}^{\{j\}}]^T\!\in\! \mathbb{C}^{(2(r-1)n+1)\times 1}$, $j=1,2,\ldots, r-1$.
\smallskip

{\bf Step 5.} Construct the $(2rn-n+1)\times (2rn-n+r+1)$ matrix $\Delta_0=[\Delta_1\;\Delta_2\;\Delta_3]$, where $\Delta_1=[T(\mathbf{c}^{\{1\}}\,;n)\;T(\mathbf{c}^{\{2\}}\,;n)\;\cdots\;
T(\mathbf{c}^{\{r-1\}}\,;n)]\in \mathbb{C}^{(2rn-n+1)\times (r-1)(n+1)}$,\;\;$\Delta_2=-T(\mathbf{b}\,; n)\in\mathbb{C}^{(2rn-n+1)\times (n+1)}$, and $\Delta_3=T(\mathbf{a}\,; (r-1)n)\in\mathbb{C}^{(2rn-n+1)\times ((r-1)n+1)}$ and then the $(2rn-n+r+1)\times (2rn-n+r+1)$ matrix $\Delta=[\Delta_1\;\Delta_2\;\Delta_3\,; \mathbf{I}\;\mathbf{0}_{r\times((r-1)n+1)}]$, where $\mathbf{I}\in\mathbb{C}^{r\times r(n+1)}$ is the $r\times r$ block identity matrix with entries $\mathbf{1}_{1\times(n+1)}$ on the block diagonal and $\mathbf{0}_{1\times(n+1)}$ elsewhere.

\smallskip

{\bf Step 6.} Perform the Cholesky factorization of the positive definite matrix $S(1)=S^+(1)\big(S^+(1)\big)^*$ and assume that $S^+(1)=[h_1\;h_2\cdots h_r]$, where $h_j\in\mathbb{C}^{r\times 1}$.

\smallskip

{\bf Step 7.} For each $j=1,2,\ldots, r$, solve the $(2rn-n+r+1)\times (2rn-n+r+1)$ system of equations
\begin{equation}\label{equDt}
\Delta X=\Lambda_j,
\end{equation}
with right-hand sides $\Lambda_j=[\mathbf{0}_{(2rn-n+1)\times 1}\;; h_j]$, and denote the respective solution  by $X_j=[x_0^{\{j\}}\,x_1^{\{j\}}\,\cdots x_{2rn-n+r}^{\{j\}}]^T$.

\smallskip

{\bf Step 8.} Set a spectral factor $S_0^+=\big(s^+_{ij}\big)_{i,j=1}^r$, where $s^+_{ij}(t)=\sum_{k=0}^n x_{(n+1)(i-1)+k}^{\{j\}} t^k$

\smallskip

{\bf Step 9.} Find $S^+$ by $S^+_0U$, where $U$ is defined by the formula \eqref{Uc}.

\smallskip

Since we know the existence of decomposition \eqref{dr7},  the solution to equation \eqref{equDt}  exists for each $j$. However it might happen that $\det\Delta=0$. Furthermore, computer simulations suggest that $\Delta$ is nonsingular whenever $\det S(t)\not= 0$ for each $t\in \mathbb{T}$ and $\Delta$ is singular whenever $\det S(t)= 0$ for some $t\in \mathbb{T}$. Therefore JLE-3 works under the additional condition $\det S(t)> 0$ for  $t\in \mathbb{T}$. If this condition holds, but zeros of $\det S$ are rather close to the boundary, the matrix $\Delta$ might become ill-conditioned. In such situations, the solutions of \eqref{equDt} are inaccurate and approximation to $S^+$ is lost. The techniques of solution of ill-conditioned systems might be useful, however we have not  investigated this question yet. As numerical simulations show in Section IX, JLE-algorithm 3 can satisfactory factorize random matrices with $r=6$ and $n=20$, which might be useful in certain applications to Mobile Communications \cite{SK2016}.

\section{Factorization of singular matrices}

Symmetric positive matrix polynomials which are chosen randomly or obtained by channel estimation in wireless communication are usually non-singular, i.e. their determinants do not vanish on $\mathbb{T}$. However, in certain optimal control and wavelet design problems, one encounters a need to factorize singular matrices. It is well known that all MSF methods have difficulties in this situation and some of them cannot handle zeros on the unit circle at all. Obviously, convergence of JLE algorithms also slows down in singular cases. However, if we fully utilize the ability of Janashia-Lagvilava's method to decompose a large scale problem into smaller parts and deal with any arising difficulties by intermediate interventions, in number of cases we can substantially improve the performance of the algorithm. In this section we demonstrate this advantage by factorizing specific singular matrices.

First, consider a test matrix from \cite{IEEE} whose spectral factorization is known beforehand:
\begin{equation}\label{IEE0}
\begin{pmatrix}2z^{-1}+6+2z&11z^{-1}+22+7z\\
7z^{-1}+22+11z&38z^{-1}+84+38z\end{pmatrix}=
\begin{pmatrix}2+z&1\\
7+5z&3+z\end{pmatrix}
\begin{pmatrix}2+z^{-1}&7+5z^{-1}\\
1&3+z^{-1}\end{pmatrix}
\end{equation}
This matrix is very simple, but its determinant, $-z^{-2}+2-z^2$,
has two double zeros on the boundary.

When data was fed into "standard" JLE-algorithm 1 with 5 iterations in scalar spectral factorization of $\det S$ by Wilson's algorithm (see Sect. 6, Procedure 1), we get 4 correct digits. When we increase the number of the iterations up to 45, the maximum optimum value, we get 7 correct digits. If we compute the determinant by the direct formula $\det S=s_{11}s_{22}-s_{12}s_{21}$, avoiding the minimal round-off errors introduced with computation of the determinant by FFT (see Section 6), then we get 14 correct digits. All these computations take less than 0.01 seconds as the matrix is very small and and it suffices to select the parameter $N$ as small as 20. We observed that Wilson's MSF algorithm (see the next section) can perform factorization \eqref{IEE0} with no more than 6 correct digits (with optimum parameter $\kappa=19$) which takes around 3 minutes.

Next we factorize a small size $2\times 2$ matrix
\begin{equation}\label{P}
S(z)=\sum_{k=-3}^3 C_k z^k=\begin{pmatrix} s_{11}(z)& s_{12}(z)\\
s_{21}(z)& s_{22}(z)\end{pmatrix},
\end{equation}
where $s_{11}(z)=-\frac{1-4\overline{\alpha}}{64}z^{-3}+\frac{1+4{\alpha}}{64}z^{-1}+1 +\frac{1+4{\alpha}}{64}z -\frac{1-4\overline{\alpha}}{64}z^{3}$;
$s_{12}(z)=\frac{\overline{\alpha}}{16}z^{-3}-\frac{{\alpha}}{16}z^{-1}+\frac{{\alpha}}{16}z -\frac{\overline{\alpha}}{16}z^{3}$; $s_{21}(z)=s_{12}(1/z)$; and
$s_{22}(z)=\frac{1-4\overline{\alpha}}{64}z^{-3}-\frac{1+4{\alpha}}{64}z^{-1}+1 -\frac{1+4{\alpha}}{64}z +\frac{1-4\overline{\alpha}}{64}z^{3}$; with $\alpha=4+\sqrt{15}$ and $\overline{\alpha}=4-\sqrt{15}$. This matrix is singular and, furthermore, its determinant has an explicit form $\det S(z)=\frac{8\overline{\alpha}-1}{4096}(z+1)^4(z-1)^4(z+i)^2(z-i)^2$. Its spectral factorization $S(z)=\sum_{k=0}^3 A_k z^k \sum_{k=0}^3 A_k^T z^{-k}$. is required for construction of the so called SA4 multiwavelet \cite{Mwt6} which possess certain nice properties. The realization of these properties depends on the accuracy by which the coefficients $A_k$ are computed.  The efforts to factorize \eqref{P} with a maximal possible accuracy by the Youla-Kazanjian method \cite{YK} is described in \cite{SA4}, where the error
$
err_1=\|S(z)-\sum\nolimits_{k=0}^3 \hat{A}_k z^k \sum\nolimits_{k=0}^3 \hat{A}_k^T z^{-k}\|=4.086\cdot 10^{-8}
$ is achieved.
(As the exact values of $A_k$ are unknown in this situation, this error is used to estimate the accuracy $\|A_k-\hat{A_k}\|$.) As we checked, this performance cannot be improved by the Wilson MSF method either. In fact, the error cannot be reduced to lower than $10^{-5}$ by the method (with optimal tuning parameter $\kappa=18$: see Section 10).

When we ran  JLE-1 with the matrix $S$ and increase the number of iterations in the scalar factorization step up to 60 (see Procedure 1), we obtain the error $err_2=4.373\cdot 10^{-5}$. However, if we cancel out the common roots in the triangular factorization \eqref{S1} and factorize the determinant $\det S$ manually we achieve the error $err_3=1.843\cdot10^{-14}$. In these computations, it is sufficient to take the tuning parameter $N=100$  and so  the consumed time is very small (less than 0.1 seconds).

In general, when a singular polynomial (with a zero on $\mathbb{T}$) is factorized in the scalar case, the best way to deal with the singularity is to factor out the zeros with unit modulus. This procedure is more demanding in the matrix case (see \cite[p. 67]{LitSpi}). The above  examples demonstrate that Janashia-Lagvilava method is capable of reducing a problem of the singularity of a spectral matrix density to the level of scalar factorization. In fact, the method has already been used to improve the coefficients of other well-known multiwavelets as well by effective factorization of related singular matrices which will be the topic of another paper.

\section{Comparison with Wilson's algorithm }

Wilson's method of MSF appeared in the 70's of the last century \cite{Wil72}, \cite{Wil78}. Since then, several authors claimed that they obtained MSF algorithms with reduced computational complexity (see \cite[p. 1077]{Kuc}, \cite[p. 206]{Kai99}). These are algorithms based on the solution of algebraic Riccati equations and some of them are implemented in Matlab. As a consequence, in our attempts to compare Janashia-Lagvilava algorithm with other existing methods of MSF, we did not originally consider the Wilson method and only concentrated our attention  on those methods which were implemented in Matlab (see \cite[Sect. VI]{IEEE}).  However, recently we learned that Prof. Rangarajan and his collaborators, who apply MSF in Neuroscience \cite{Dhamala}, \cite{Ranga}, developed an  efficient implementation of Wilson's method which works rather fast.

This implementation takes data matrix in frequency domain. Nevertheless, this idea can be easily translated for matrices given in time domain. In particular, for a matrix \eqref{polm} with given coefficients $C_k$, $k=0,1,\ldots,N$, we select $\kappa$ as a tuning parameter and find $2^\kappa$ values of the matrix function $S$ in DFT nodes: $S(t_0),\dots,S(t_{2^\kappa})$, where $t_j=\exp\left(\frac{2\pi ij}{2^\kappa}\right)$. Then we use the Wilson's recurrent formula
\begin{equation}\label{Wil}
S^+_{k+1}=S^+_k\left[(S^+_k)^{-1}S(S^+_k)^{-*}+I\right]^+
\end{equation}
with initial data $S_0=\sqrt{C_0}$. After performing sufficient iterations, we return back to the time domain and approximately compute  the coefficients $A_k$ of \eqref{polm1}. Here, like other minor improvements we introduced in the implementation of Wilson's method, we empirically observed that the upper triangular constant matrix $S_\tau$ in formula (3.2) in \cite{Wil78} can be omitted in \eqref{Wil}.
Such implementation of Wilson's algorithm essentially works as efficient as JLE-1 and frequently better than JLE-2. In addition, a flexible combination of Janashia-Lagvilava and Wilson methods can be sometimes useful.

\section{Numerical simulations}

The computer code for implementation of JLE-algorithms was written in Matlab in order to test them numerically.  A laptop
with characteristics Intel(R) Core(TM) i7-4600U CPU (2 cores, 4 threads), 2.40GHz, RAM 8.00Gb was used and some of the tests were performed on the HPC cluster ``Dalma" at NYUAD.

For all numerical simulations of MSF algorithms randomly selected polynomial matrices have been used. Namely, for  given matrix dimension $r$ and polynomial degree $n$, a random polynomial matrix $\sum_{k=0}^n A_kt^k$, $A_k\in [-1,1]^{r\times r}$, has been chosen, and positive definite (on $\mathbb{T}$) matrix polynomial $S(t)=\sum_{k=0}^n A_kt^k\sum_{k=0}^n A_k^*t^{-k}$ has been approximately factorized. In rare occasions, which are emphasized below, some deterministic efforts have been introduced in order to artificially improve the properties of $S$. The error
\begin{equation}\label{err}
err=\|S-\hat{S}^+(\hat{S}^+)^*\|
\end{equation}
is used to estimate the accuracy of the factorization since there is no other way to decide how close is $\hat{S}^+$ ro $S^+$.

The basic problem in order to demonstrate the most effective performance of the constructed algorithms was an empirical selection of tuning parameters ($N$ for JLE-1, $N$ and $\kappa$ for JLE-2, and $\kappa$ and the number of iterations for Wilson's algorithm) which would make an optimal trade-off between the available memory, the computation time and the accuracy.

For realistic applications, automatic selection of the optimal tuning parameters during the factorization remains a challenging problem.

When different algorithms are compared, it is assumed that they were run with the same data.

We start with JLE-3 which has the advantage that it contains no tuning parameters. Below we demonstrate its performance within the range of polynomial matrices for which it is applicable. The tuning parameters in JLE-1 and Wilson have been selected so as to achieve the same accuracy as in JLE-3. Beyond the indicated range of matrix dimension $m$ and polynomial degree $n$ the accuracy \eqref{err} of JLE-3 becomes unsatisfactory. (In all tables below, $r\times n$ indicates that a $r\times r$ test matrix was selected with Laurent polynomial entries of degree $n$ having nonzero coefficients indexed from $-n$ to $n$).

\begin{center}
{\footnotesize Table I\\Performance of JLE-3\par}
\end{center}
{\fontsize{8}{8pt}\selectfont
\begin{center}
\begin{tabular}{|c|c|c|c|c|c|c|c|}
\hline &&&&&&\\[-2mm]
&matrix&time&accu-&matr.&time&accu-\\
&size&$sec$&racy&size&$sec$&racy\\
\hline &&&&&&\\[-2mm]
JLE-3&$4\times 30$&0.\,052&$10^{-8}$&$6\times 20$&0.\,051&$10^{-6}$\\
\hline &&&&&&\\[-2mm]
JLE-1&--&0.\,315&$10^{-8}$&--&0.\,576&$10^{-6}$\\
\hline &&&&&&\\[-2mm]
Wilson&--&1.\,108&$10^{-8}$&--&0.\,694&$10^{-6}$\\
\hline &&&&&&\\[-2mm]
JLE-3&$8\times 10$&0.\,051&$10^{-6}$&$10\times 5$&0.\,044&$10^{-6}$\\
\hline &&&&&&\\[-2mm]
JLE-1&--&0.\,359&$10^{-6}$&--&0.\,352&$10^{-6}$\\
\hline &&&&&&\\[-2mm]
Wilson&--&0.\,419&$10^{-6}$&--&0.\,360&$10^{-6}$\\
\hline
\end{tabular}
\end{center}}

\smallskip

Next we compare JLE-1 and Wilson within the range of matrices where JLE-1 operates well. The tuning parameter $N=5mn$ has been taken for $m$th recursion in JLE-1 and $\kappa$ has been selected in Wilson so as to achieve the same accuracy as in JLE-1.
\begin{center}
{\footnotesize Table II\\Comparision of JLE-1 and Wilson\par}
\end{center}
{\fontsize{8}{8pt}\selectfont
\begin{center}
\begin{tabular}{|c|c|c|c|c|c|}
\hline &&&&\\[-2mm]
&matrix&tuning&time&accuracy\\
&size&parameters&$sec$&\\
\hline &&&&\\[-2mm]
JLE-1&$10\times 100$&N = 500 m&6.35&$1.63\cdot10^{-7}$\\
\hline &&&&\\[-2mm]
Wilson&--&$\kappa = 12$; Iter = 23&$7.91$&$1.42\cdot10^{-7}$\\
\hline &&&&\\[-2mm]
JLE-1&$15\times 20$&N = 100 m&2.67&$6.12\cdot10^{-8}$\\
\hline &&&&\\[-2mm]
Wilson&--&$\kappa = 11$; Iter = 25 &3.95&$1.76 \cdot10^{-8}$\\
\hline
\end{tabular}
\end{center}}

\smallskip

Next we factorize random $100\times 100$ matrices (with  polynomial degree $n=30$) by JLE-2 and Wilson. We tried to factorize such matrices with accuracy that is acceptable in practice, namely  $error=10^{-4}$,  and  selected the tuning parameters accordingly. A substantial drop in the accuracy has been observed at the final step of recursion  $m=100$ in JLE-2 and it was observed that Wilson can factorize the $99\times 99$ leading submatrix of $S$ much more easily than $S$ itself. We  empirically explain this phenomenon  by the following reason: the probability for zeros of $\det S_{[m]}$ to be very close to $\mathbb{T}$ (in which case all spectral factorization algorithms become slowly convergent) is higher for $m=r$ than for $m<r$  (however no theoretical proofs has been attempted). Therefore, in a variant of our implementation, we have combined JLE-2 by Wilson  which resulted in certain improvements.
\begin{center}
{\footnotesize Table III\\Comparision of JLE-2 and Wilson\par}
\end{center}
{\fontsize{8}{8pt}\selectfont
\begin{center}
\begin{tabular}{|c|c|c|c|c|}
\hline &&&\\[-2mm]
$100\times 30$&tuning parameters&time&accuracy\\[+1mm]
\hline &&&\\[-2mm]
JLE-2&\tiny{\!$N=2400e^{0.15(m-100)}+80e^{0.01m}\! $}&94.1&$5\cdot10^{-4}$\\
\hline &&&\\[-2mm]
Wilson&$\kappa = 11$; Iter = 14&$85.3$&$1\cdot10^{-4}$\\
\hline &&&\\[-2mm]
Wil+JLE&{$99\times99$ by Wilson+JLE-2} &43.3&$2\cdot10^{-4}$\\
\hline
\end{tabular}
\end{center}}

\smallskip

 When we added artificially $I_r$  to a random matrix $S$ in order to avoid zeros close to $\mathbb{T}$, we achieved the same accuracy within improved computation time. We display the results below.

\begin{center}
{\footnotesize Table IV\\Comparision of JLE-2 and Wilson\par}
\end{center}
{\fontsize{8}{8pt}\selectfont
\begin{center}
\begin{tabular}{|c|c|c|c|c|}
\hline &&&\\[-2mm]
$S_{rand}+I\!\!$&tuning parameters&time&accuracy\\[+1mm]
\hline &&&\\[-2mm]
JLE-2&\tiny{\!$N\approx 100e^{0.02m} $}&24.1&$4\cdot10^{-4}$\\
\hline &&&\\[-2mm]
Wilson&$\kappa = 8$; Iter = 9&$6.57$&$3\cdot10^{-4}$\\
\hline &&&\\[-2mm]
Wil+JLE&{$99\times99$ by Wilson+JLE-2} &6.23&$5\cdot10^{-4}$\\
\hline
\end{tabular}
\end{center}}

\smallskip

In the end we demonstrate that ``good" matrices of dimension as large as $700\times 700$ can be factorized with accuracy $error=10^{-3}$ which is acceptable
 in practice and within the available computer memory (120GB of one node at ``Dalma" in our situation). With respect to time usage, the advantage of Wilson's MSF method is evident in this case. The reason is that JLE-2 requires the tuning parameter $N$ to be selected  very large at the last recursive  steps in order to achieve the given accuracy. However, JLE-2 algorithm still can be invoked to analyze and overcome the problem when Wilson's method is unable to factorize a matrix obtained from real applications.

 \begin{center}
{\footnotesize Table V\\Factorization of large matrices\par}
\end{center}
{\fontsize{8}{8pt}\selectfont
\begin{center}
\begin{tabular}{|c|c|c|c|c|c|}
\hline &&&&\\[-2mm]
$700\times 30$&tuning &time &accu-&RAM\\
$S_{rnd}+I$&parameters&(hours)&racy&\\[+1mm]
\hline &&&&\\[-2mm]
JLE-2&\tiny{\!$N\approx 200e^{0.003m}\! $}&$2:50$&$4\cdot10^{-3}$&80GB\\
\hline &&&&\\[-2mm]
Wilson&$\kappa = 11$; Iter = 17&$1:08$&$3\cdot10^{-3}$&120GB\\
\hline
\end{tabular}
\end{center}}

\section{Conclusions}

Matrix spectral factorization is widely used in modern control theory and wireless communications. Furthermore, improved algorithms of MSF may lead to new  areas to which they could be successfully applied. In the present paper, we consider three different algorithms based on Janashia-Lagvilava method, which may be competitive with other existing MSF algorithms. A general description of their computational capabilities, as well as a comparison to Wilson's MSF algorithm, are provided by means of numerical simulations.

\section{Acknowledgments}
The authors are thankful for an opportunity to run part of the tests using the High Performance Computing resources at New York University Abu Dhabi.

\def\cprime{$'$}
\providecommand{\bysame}{\leavevmode\hbox to3em{\hrulefill}\thinspace}
\providecommand{\MR}{\relax\ifhmode\unskip\space\fi MR }
\providecommand{\MRhref}[2]{%
  \href{http://www.ams.org/mathscinet-getitem?mr=#1}{#2}
}
\providecommand{\href}[2]{#2}


\begin{thebibliography}{10}

\bibitem{Bott13}
A.~B{\"o}ttcher and M.~Halwass, \emph{A {N}ewton method for canonical
  {W}iener-{H}opf and spectral factorization of matrix polynomials}, Electron.
  J. Linear Algebra \textbf{26} (2013), 873--897. \MR{3192406}

\bibitem{Dhamala}
M.~Dhamala, G.~Rangarajan, and M~Ding, \emph{Analyzing information flow in
  brain networks with nonparametric granger causality}, NeuroImage \textbf{41}
  (2008), 354–--362.

\bibitem{Ranga}
\bysame, \emph{Estimating granger causality from fourier and wavelet transforms
  of time series data}, Physical Review Letters \textbf{100} (2008), 018701.

\bibitem{E}
L.~Ephremidze, \emph{An elementary proof of the polynomial matrix spectral
  factorization theorem}, Proc. Roy. Soc. Edinburgh Sect. A \textbf{144}
  (2014), no.~4, 747--751. \MR{3233753}

\bibitem{EJL11}
L.~Ephremidze, G.~Janashia, and E.~Lagvilava, \emph{On approximate spectral
  factorization of matrix functions}, J. Fourier Anal. Appl. \textbf{17}
  (2011), no.~5, 976--990. \MR{2838115 (2012h:47039)}

\bibitem{EL10}
L.~Ephremidze and E.~Lagvilava, \emph{Remark on outer analytic
  matrix-functions}, Proc. A. Razmadze Math. Inst. \textbf{152} (2010), 29--32.
  \MR{2663529}

\bibitem{EL2014}
\bysame, \emph{On compact wavelet matrices of rank {$m$} and of order and
  degree {$N$}}, J. Fourier Anal. Appl. \textbf{20} (2014), no.~2, 401--420.
  \MR{3200928}

\bibitem{ESS}
L.~Ephremidze, N.~Salia, and I.~Spitkovsky, \emph{Some aspects of a novel
  matrix spectral factorization algorithm}, Proc. A. Razmadze Math. Inst.
  \textbf{166} (2014), 49--60. \MR{3300615}

\bibitem{SFLP}
T.~N.~T. Goodman, Ch.~A. Micchelli, G.~Rodriguez, and S.~Seatzu, \emph{Spectral
  factorization of {L}aurent polynomials}, Adv. Comput. Math. \textbf{7}
  (1997), no.~4, 429--454. \MR{1470294}

\bibitem{HelLow58}
H.~Helson and D.~Lowdenslager, \emph{Prediction theory and {F}ourier series in
  several variables}, Acta Math. \textbf{99} (1958), 165--202. \MR{0097688 (20
  \#4155)}

\bibitem{Jaf}
A.~Jafarian and J.~G. McWhirter, \emph{A novel method for multichannel spectral
  factorization}, Proc. Europ. Signal Process. Conf. (2012), 27--31.

\bibitem{JL99}
G.~Janashia and E.~Lagvilava, \emph{A method of approximate factorization of
  positive definite matrix functions}, Studia Math. \textbf{137} (1999), no.~1,
  93--100. \MR{1735630 (2000m:15015)}

\bibitem{IEEE}
G.~Janashia, E.~Lagvilava, and L.~Ephremidze, \emph{A new method of matrix
  spectral factorization}, IEEE Trans. Inform. Theory \textbf{57} (2011),
  no.~4, 2318--2326. \MR{2809092 (2012d:65077)}

\bibitem{Kai99}
T.~Kailath, B.~Hassibi, and A.~H. Sayed, \emph{Linear estimation},
  Prentice-Hall, Inc., Englewood Cliffs, N.J., 1999, Prentice-Hall Information
  and System Sciences Series.

\bibitem{SA4}
V.~Kolev, T.~V. Cooklev, and F.~Keinert, \emph{Matrix spectral factorization -
  {SA4} multiwavelet}, Preprint.

\bibitem{Kuc}
V.~Ku{\v{c}}era, \emph{Factorization of rational spectral matrices: A survey of
  methods}, in Proc. IEEE Int. Conf. Control, Edinburgh \textbf{2} (1991),
  1074--1078.

\bibitem{LitSpi}
G.~S. Litvinchuk and I.~M. Spitkovskii, \emph{Factorization of measurable
  matrix functions}, Operator Theory: Advances and Applications, vol.~25,
  Birkh\"auser Verlag, Basel, 1987, Translated from the Russian by Bernd
  Luderer, With a foreword by Bernd Silbermann. \MR{1015716}

\bibitem{SayKai}
A.~H. Sayed and T.~Kailath, \emph{A survey of spectral factorization methods},
  Numer. Linear Algebra Appl. \textbf{8} (2001), no.~6-7, 467--496, Numerical
  linear algebra techniques for control and signal processing. \MR{1848590
  (2002j:93039)}

\bibitem{SK2016}
A.~B. Sergienko and V.~P Klimentyev, \emph{Scma detection with channel
  estimation error and resource block diversity}, in Proc. Int. Siberian Conf.
  Control and Communications, (SIBCON) (2016), DOI:
  10.1109/SIBCON.2016.7491765.

\bibitem{Mwt6}
J.~Y. Tham, L.~Shen, S.~L. Lee, and H.~H. Tan, \emph{A general approach for
  analysis and application of discrete multiwavelet transforms}, IEEE Trans.
  Signal Process. \textbf{48} (2000), no.~2, 457--464. \MR{1746064}

\bibitem{RoyalA}
X.~Wen, G.~Rangarajan, and M.~Ding, \emph{Multivariate granger causality: an
  estimation framework based on factorization of the spectral density matrix},
  Phil. Trans. R. Soc. A 371: 20110610. (2013).

\bibitem{Wie57}
N.~Wiener and P.~Masani, \emph{The prediction theory of multivariate stochastic
  processes. {I}. {T}he regularity condition}, Acta Math. \textbf{98} (1957),
  111--150. \MR{0097856 (20 \#4323)}

\bibitem{Wie58}
\bysame, \emph{The prediction theory of multivariate stochastic processes.
  {II}. {T}he linear predictor}, Acta Math. \textbf{99} (1958), 93--137.
  \MR{0097859 (20 \#4325)}

\bibitem{Wil1}
G.~Wilson, \emph{Factorization of the covariance generating function of a pure
  moving average process}, SIAM J. Numer. Anal. \textbf{6} (1969), 1--7.
  \MR{0253561}

\bibitem{Wil72}
G.~Tunnicliffe Wilson, \emph{The factorization of matricial spectral
  densities}, SIAM J. Appl. Math. \textbf{23} (1972), 420--426. \MR{0331843}

\bibitem{Wil78}
\bysame, \emph{A convergence theorem for spectral factorization}, J.
  Multivariate Anal. \textbf{8} (1978), no.~2, 222--232. \MR{497596}

\bibitem{YK}
D.~C. Youla and N.~N. Kazanjian, \emph{Bauer-type factorization of positive
  matrices and the theory of matrix polynomials orthogonal on the unit circle},
  IEEE Trans. Circuits and Systems \textbf{CAS-25} (1978), no.~2, 57--69.
  \MR{0469461}

\end{thebibliography}
\end{document}